\documentclass{article}
\usepackage{graphicx} 
\usepackage[T1]{fontenc}
\usepackage[latin9]{inputenc}
\usepackage{float}
\usepackage{amsthm}
\usepackage{amsmath}
\usepackage{mathrsfs}
\usepackage{amssymb}
\usepackage{graphicx}
\usepackage[]{wrapfig}
\usepackage[export]{adjustbox}
\usepackage{subfigure}
\usepackage{enumitem}
\usepackage[english]{babel}
\usepackage{mathtools}
\usepackage{newlfont}
\usepackage{indentfirst}
\usepackage{booktabs}
\usepackage{multirow}
\usepackage{lineno}

\usepackage{hyperref}
\usepackage[capitalise]{cleveref}

\theoremstyle{plain} 
\newtheorem{thm}{Theorem}[section] 
\newtheorem{cor}[thm]{Corollary}

\newtheorem{prob}[thm]{Problem}
 
\newtheorem{conj}[thm]{Conjecture}

\newtheorem{obs}[thm]{Observation}

\usepackage{xcolor}

\hypersetup{colorlinks = true, linkcolor = black, citecolor = black, urlcolor = blue}

\title{The Gray graph is pseudo 2-factor isomorphic }
\date{}
\author{
Mari\'en Abreu \thanks{Department of Basic and Applied Sciences, Universit\`a degli Studi della Basilicata, 85100 Potenza, Italy. 
\protect\href{mailto:marien.abreu@unibas.it}{\protect\nolinkurl{marien.abreu@unibas.it}}}
\and Jan Goedgebeur \thanks{Department of Computer Science, KU Leuven Campus Kulak-Kortrijk, 8500 Kortrijk, Belgium. 
\protect\href{mailto:jan.goedgebeur@kuleuven.be}{\protect\nolinkurl{jan.goedgebeur@kuleuven.be}},
\protect\href{mailto:jorik.jooken@kuleuven.be}{\protect\nolinkurl{jorik.jooken@kuleuven.be}},
 and \protect\href{mailto:tibo.vandeneede@kuleuven.be}{\protect\nolinkurl{tibo.vandeneede@kuleuven.be}}} 
 \thanks{Department of Mathematics, Computer Science and Statistics, Ghent University, 9000 Ghent, Belgium.}
\and Jorik Jooken
 \footnotemark[2] 
\and Federico Romaniello \thanks{Department for Humanistic, Scientific and Social Innovation, Universit\`a degli Studi della Basilicata, 85100 Potenza, Italy. 
\protect\href{mailto:federico.romaniello@unibas.it}{\protect\nolinkurl{federico.romaniello@unibas.it}}}  
\and Tibo Van den Eede
 \footnotemark[2]
}

\begin{document}

\maketitle

\begin{abstract}

A graph is \textit{pseudo 2-factor isomorphic} if all of its 2-factors have the same parity of number of cycles. 

Abreu et al.\ [J.\ Comb.\ Theory, Ser.\ B.\ 98 (2008) 432--442] conjectured that $K_{3,3}$, the Heawood graph and the Pappus graph are the only essentially 4-edge-connected pseudo 2-factor isomorphic cubic bipartite graphs. 

This conjecture was disproved by Goedgebeur [Discr.\
Appl.\ Math.\ 193 (2015) 57--60] who constructed a counterexample $\mathcal{G}$ (of girth 6) on 30 vertices. Using a computer search, he also showed that this is the only counterexample up to at least 40 vertices and that there are no counterexamples of girth greater than 6 up to at least 48 vertices.

In this manuscript, we show that the Gray graph -- which has 54 vertices and girth 8 -- is also a counterexample to the pseudo 2-factor isomorphic graph conjecture. Next to the graph $\mathcal{G}$, this is the only other known counterexample.  Using a computer search, we show that there are no smaller counterexamples of girth 8 and show that there are no other counterexamples up to at least 42 vertices of any girth. 

Moreover, we also verified that there are no further counterexamples among the known censuses of symmetrical graphs. 

Recall that a graph is \textit{2-factor Hamiltonian} if all of its 2-factors are Hamiltonian cycles. 
As a by-product of the computer searches performed for this paper, we have verified that the \textit{$2$-factor Hamiltonian conjecture} of Funk et al.\ [J.\ Comb.\ Theory, Ser.\ B.\ 87(1) (2003) 138--144], which is still open, holds for cubic bipartite graphs of girth at least 8 up to 52 vertices, and up to 42 vertices for any girth. 

\bigskip


\end{abstract}

\section{Introduction}
\begin{figure}[h!]
    \centering
    \includegraphics[width=1\textwidth]{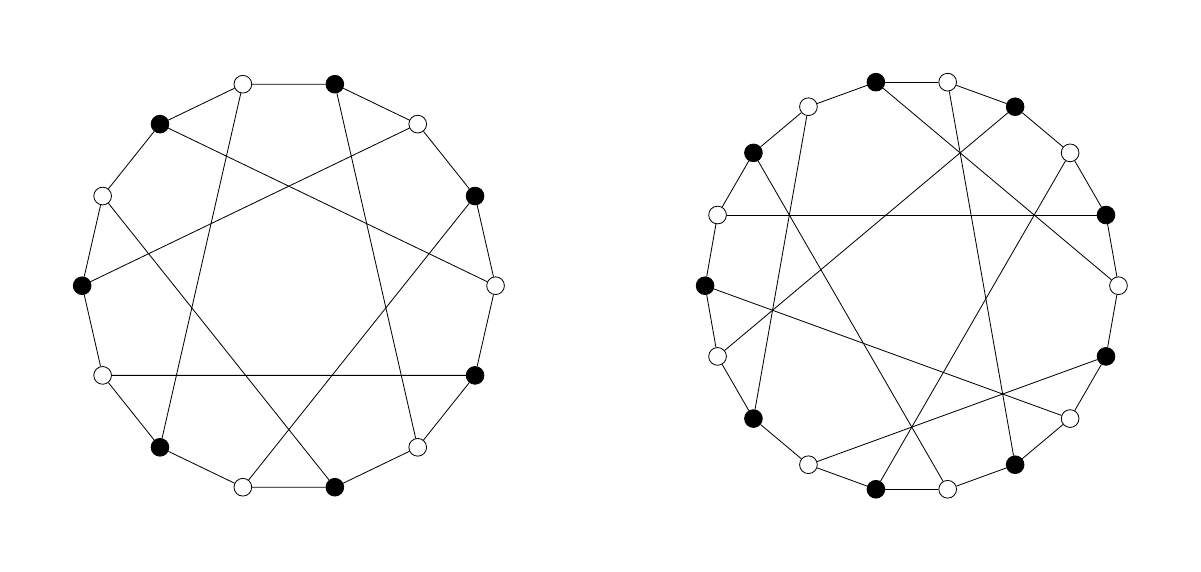} 
    \caption{The Heawood graph on the left and the Pappus graph on the right.}
    \label{fig:heawpapp}
\end{figure}

A spanning cycle in a graph is a \textit{Hamiltonian cycle}, 
where \textit{spanning} means that it contains all vertices of the graph.
The problem of finding a Hamiltonian cycle is well known to be NP-complete~\cite{Kar72}. The study of \textit{Hamiltonicity} in graphs, i.e., whether or not a graph contains a Hamiltonian cycle, in particular for regular or cubic graphs, has been considered by many authors with beautiful results and conjectures that can be found in the literature (see e.g.,~\cite{RG91,RG14}). Since a good characterisation, say with a clear structure, coming from a (finite and not too long) set of properties, of Hamiltonian graphs is not likely to exist, there are some related problems for which it seems more feasible to understand the nature of the objects involved. Let us recall that a \emph{$k$-factor} of $G$ is defined as a $k$-regular spanning subgraph of $G$ (not necessarily connected). Thus, a \emph{$2$-factor} of a graph $G$ is a $2$-regular spanning subgraph of $G$, i.e.\ it is a Hamiltonian cycle or a set of disjoint cycles covering all vertices of the graph. A graph $G$ is \emph{$2$-factor Hamiltonian (isomorphic)} if all of its $2$-factors are Hamiltonian cycles (resp.\ isomorphic). Examples of $2$-factor Hamiltonian graphs are the complete graphs $K_4$ and $K_5$; the complete bipartite graph $K_{3,3}$; the Heawood graph $H_0$, i.e.\ the Levi (incidence) graph of the Fano Plane $PG(2,2)$ which is also the $7_3$ symmetric configuration (cf.\ Figure~\ref{fig:heawpapp} left). The latter two examples, which happen to be bipartite, are quite remarkable, as we will see briefly. Recall that a graph is \textit{bipartite} if its vertex set can be partitioned into two sets in such a way that each edge has one end in each of them. 
An example of a $2$-factor isomorphic graph which is neither Hamiltonian nor bipartite is the Petersen graph.
In the bipartite case, $2$-factor Hamiltonian graphs have been studied and partially characterised in~\cite{FJLS}, where the following conjecture was stated:

\begin{conj}\cite[Conjecture $3.2$]{FJLS}\label{ConjFJLS}
Let $G$ be a $2$-factor Hamiltonian $k$-regular bipartite graph. Then either $k=2$ and $G$ is a cycle or $k=3$ and $G$ can be obtained from $K_{3,3}$ and $H_0$ by repeated star products.
\end{conj}

A graph $G$ is a \emph{star product} of the graphs $G_1$ and $G_2$ if and only if there is an $x \in V(G_1)$ with neighbours $x_1, x_2, x_3$ in $G_1$ and a $y \in V(G_2)$ with neighbours $y_1, y_2, y_3$ in $G_2$ such that $G = (G_1 - x) \cup (G_2 - y) \cup \{(x_1, y_1),(x_2, y_2),(x_3, y_3)\}$. 

In~\cite{FJLS} it was also proved that there are no $2$-factor Hamiltonian $k$-regular bipartite graphs for $k \ge 4$, and that for $k=3$ these graphs are $3$-connected and have a number of vertices congruent to 2 modulo 4.

In a cubic graph, the three edges incident with a vertex constitute a $3$-edge-cut because their removal leaves an isolated vertex, and is called \emph{trivial}, other edge-cuts being \emph{non-trivial}. 
A set $S$ of edges of a graph $G$ is said to be a \emph{cyclic edge-cut} if removing $S$ from $G$ results in two components, each containing a cycle. 
A graph $G$ is said to be \emph{cyclically $m$-edge-connected} if every cyclic edge-cut in $G$ has size at least $m$.
A cubic graph is said to be \emph{essentially $4$-edge-connected} if it is $3$-connected and every $3$-edge-cut is trivial.
For an integer $k$ with $k \leq 4$, a cubic graph is essentially $k$-edge-connected if and only if it is cyclically $k$-edge-connected (see e.g.~\cite{CCOVY}). 
Please note that graphs obtained from $K_{3,3}$ and $H_0$ by repeated star products are not cyclically 4-edge-connected because of the $3$-edge-cut produced by the definition of the star product.

The star product of graphs has an inverse operation called  \textit{3-cut reduction} which consists of removing a non-trivial $3$-edge-cut and adding a vertex of degree 3 to each of the two remaining components. By repeatedly performing this operation on a cubic graph $G$, we arrive at the essentially $4$-edge-connected \textit{constituents} of $G$.

The results in~\cite{DL01,DL02} about \textit{minimally 1-factorable graphs}, imply that a counterexample to Conjecture~\ref{ConjFJLS} is cyclically $4$-edge-connected and has girth at least $6$. 
So to prove Conjecture \ref{ConjFJLS}, it would be sufficient to show:

\begin{conj}[\cite{FJLS}]\label{ConjHeaOnly}
The Heawood graph is the only $2$-factor Hamiltonian cyclically $4$-edge-connected cubic bipartite graph of girth at least $6$.
\end{conj}

Recently in~\cite{GJW} it was proved that Conjecture~\ref{ConjFJLS} is also equivalent to \cite[Conjecture 1.8]{GJW}:
\textit{The Heawood graph and $K_{3,3}$ are the only $2$-factor Hamiltonian, cubic braces.} 
Moreover, it can be narrowed to \cite[Conjecture 1.9]{GJW}: \textit{All non-Pfaffian, cubic braces of girth at least $6$, of order
congruent to $2$ modulo $4$, are not $2$-factor Hamiltonian.}  For definitions of braces and Pfaffian, please refer to~\cite{GJW}.

Another equivalence of Conjecture \ref{ConjFJLS} appears in the context of extending perfect matchings to Hamiltonian cycles and is stated in \cite[Conjecture 7]{RZ} as: \textit{Every bipartite cyclically 4-edge-connected cubic Perfect Matching Hamiltonian graph with girth at least 6, except the Heawood graph, admits a perfect matching which can be extended to a Hamiltonian cycle in exactly one way.} For definitions of \textit{perfect matching} and \textit{Perfect Matching Hamiltonian graph} please refer to \cite{RZ}.

In~\cite{AlFJLS} it was conjectured that in the cubic bipartite case, $2$-factor isomorphic graphs are also $2$-factor Hamiltonian, but this turned out to be false for the $2$-edge connected case, with counterexamples presented in~\cite{ADJLS}. In the latter paper the previous results about $2$-factor isomorphic graphs and $2$-factor Hamiltonian graphs are extended to the more general family of pseudo $2$-factor isomorphic graphs, which is the main topic of our study in this paper. A graph is \textit{pseudo $2$-factor isomorphic} if all of its $2$-factors have the same parity of number of cycles. In particular, in~\cite{ADJLS} a partial characterisation is given and it is conjectured that:

\begin{conj}\cite[Conjecture $3.5$]{ADJLS}\label{cong1}
Let $G$ be a $3$-edge-connected cubic bipartite graph. Then $G$ is pseudo 
$2$-factor isomorphic if and only if $G$ can be obtained from $K_{3,3}$, the Heawood graph or the Pappus graph by repeated star products.
\end{conj}

Similarly to the Heawood graph, which arises as the Levi graph of the $7_3$ Fano configuration, the Pappus graph is the Levi graph of the $9_3$ Pappus configuration (cf.\ Figure~\ref{fig:heawpapp} right). They are both cubic bipartite graphs of girth 6. Each bipartite graph of girth at least 6 is the Levi graph of some abstract incidence structure, in particular, if such a graph has order $\nu$ and is $k$-regular, it is the Levi graph of a  $\nu_k$ symmetric configuration. Please recall that a symmetric configuration $\nu_k$ is an incidence structure with $\nu$ points and lines, each point being incident with $k$ lines and each line containing $k$ points such that two points lie in at most one line. This last condition prevents the Levi graph of the symmetric configuration from containing 4-cycles.
Hence, by the results in~\cite{ALS}, a counterexample to Conjecture~\ref{cong1} must arise as the incidence graph of a \textit{reducible} symmetric configuration $\nu_3$. For definitions of \textit{reducible symmetric configuration}, please refer to~\cite{ALS}.

The \textit{type} of a $2$-factor in a graph is the tuple listing the lengths of its cycles in increasing order. So, for example the unique $2$-factor types of $K_4$, $K_5$, $K_{3,3}$, $H_0$ and the Petersen graph are respectively $(4), \, (5), \, (6), \, (14)$ and $(5,5)$; while the Pappus graph has $2$-factors of two types, precisely $(18)$ and $(6,6,6)$.

The classes of graphs mentioned so far have also led to the study of other conditions on the $2$-factors of a graph, such as having only odd (even) cycles in their $2$-factors, called \textit{odd (even) $2$-factored} graphs. Examples of odd $2$-factored graphs of course start with the Petersen graph, with $2$-factors only of type $(5,5)$ and which is a snark. In fact, further examples of odd $2$-factored graphs can be found among snarks (cf.~\cite{ALRS,LR}).  While bipartite graphs are always even $2$-factored, there are also non-bipartite ones which appear in~\cite{AGLRZ} with applications to extending perfect matchings to Hamiltonian cycles in graphs.

\

Similarly to the $2$-factor Hamiltonian case, it follows that Conjecture \ref{cong1} holds only if Conjecture \ref{cong2} (below) holds because graphs obtained by repeated star products are not essentially 4-edge-connected, i.e.\ by repeated 3-cut reductions we arrive at the constituents which are essentially 4-edge connected.

\begin{conj}\cite[Conjecture $3.6$]{ADJLS}\label{cong2}
Let $G$ be an essentially $4$-edge-connected pseudo $2$-factor isomorphic cubic bipartite graph. Then $G$ must be $K_{3,3}$, the Heawood graph or the Pappus graph.
\end{conj}
\begin{figure}[h!]
    \centering
    \includegraphics[width=.4\textwidth]{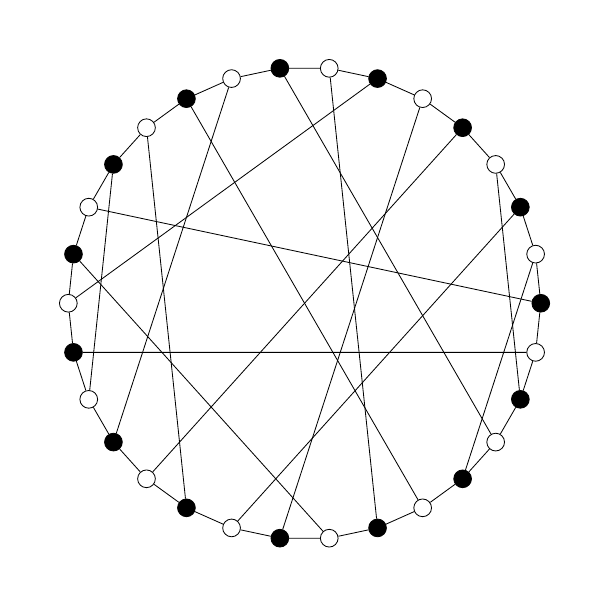} 
    \caption{The pseudo 2-factor isomorphic graph $\mathcal{G}$ on 30 vertices.}
    \label{fig:jan30}
\end{figure}

Moreover, by \cite[Theorem 3.15]{ADJLS} in a 3-edge-connected pseudo 2-factor isomorphic bipartite graph $G$, any  4-cycle $C$ is contained in a constituent of $G$ which is isomorphic to $K_{3,3}$. Therefore, for the rest of the paper, our search will be restricted to graphs of girth at least 6.

\begin{figure}[h!]
    \centering
    \includegraphics[width=.8\textwidth]{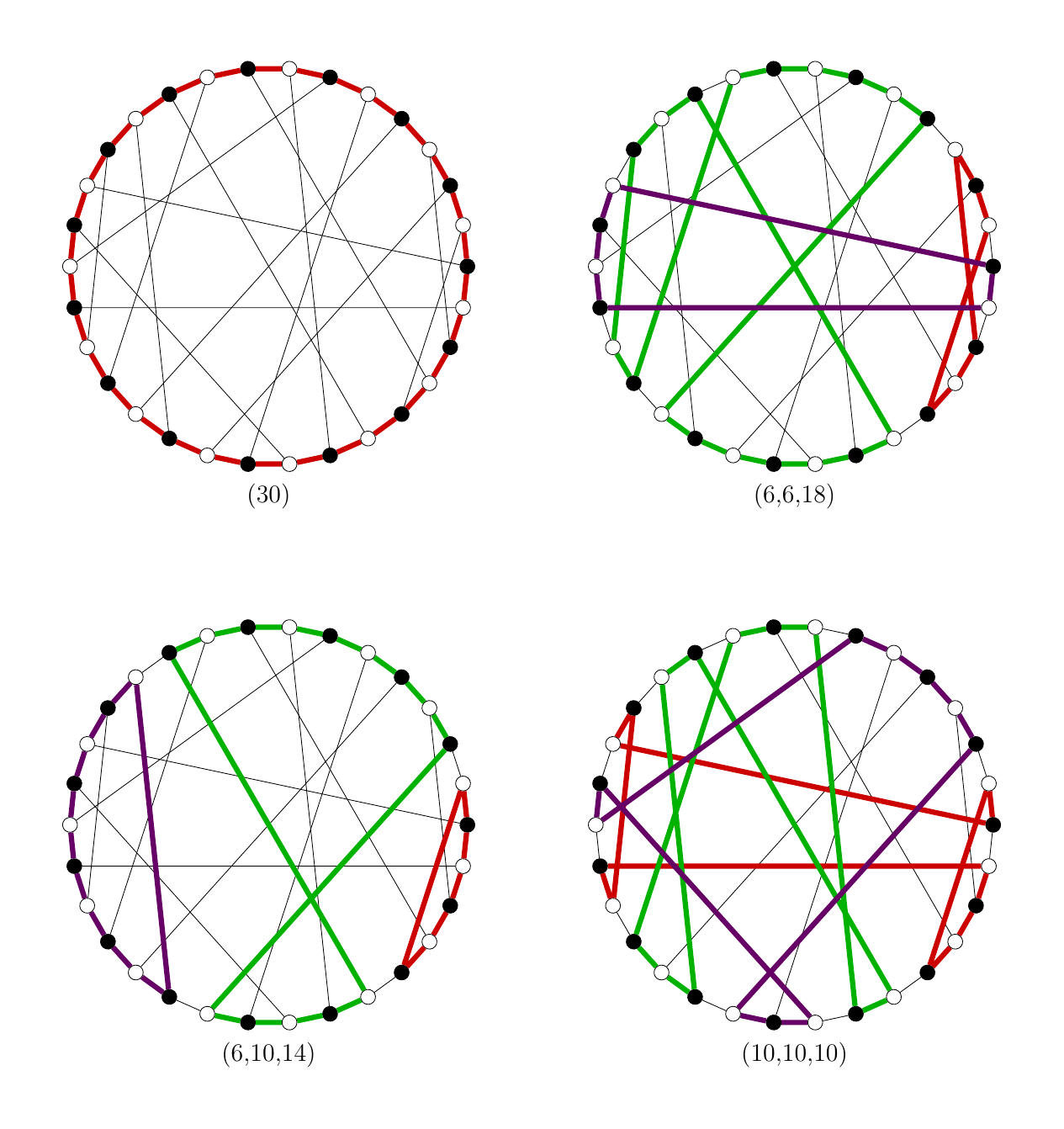} 
    \caption{The counterexample $\mathcal{G}$ on 30 vertices, with its 2--factor types highlighted.}
    \label{fig:jan}
\end{figure}

In~\cite{JG} the second author refuted Conjecture~\ref{cong2} (and consequently also Conjecture~\ref{cong1}) by constructing a counterexample using a computer search. In particular, he generated all cubic bipartite graphs with girth at least 6 up to 40 vertices and all cubic bipartite graphs with girth at least 8 up to 48 vertices and tested which of those graphs is pseudo 2-factor isomorphic. This exhaustive search yielded a \emph{new} pseudo 2-factor isomorphic graph $\mathcal{G}$ on $30$ vertices with girth 6, see Figure~\ref{fig:jan30}. $\mathcal{G}$ is essentially $4$-edge-connected, has cyclic edge-connectivity 6, automorphism group isomorphic to $(\mathbb{Z}_3 \times \mathbb{Z}_3)\rtimes(D_4 \times \mathbb{Z}_2)$ of size 144, is neither vertex-transitive nor edge-transitive, has 312 2-factors and the types of its 2-factors are: $(6,6,18)$, $(6,10,14)$, $(10,10,10)$ and $(30)$, see Figure~\ref{fig:jan}. Thus, $\mathcal{G}$ is a counterexample for the above Conjecture~\ref{cong2} (but not for Conjecture~\ref{ConjHeaOnly} as it is not 2-factor Hamiltonian), and it is even $2$-factored, being bipartite.\\

Recently, in~\cite{AFLR} a geometric construction of $\mathcal{G}$ has been obtained from the Heawood graph and the generalised Petersen graph $GP(8,3)$, respectively Levi graphs of the Fano $7_3$ and the M\"obius-Kantor $8_3$ symmetric configurations.
\begin{figure}[h!]
    \centering
    \includegraphics[width=.6\textwidth]{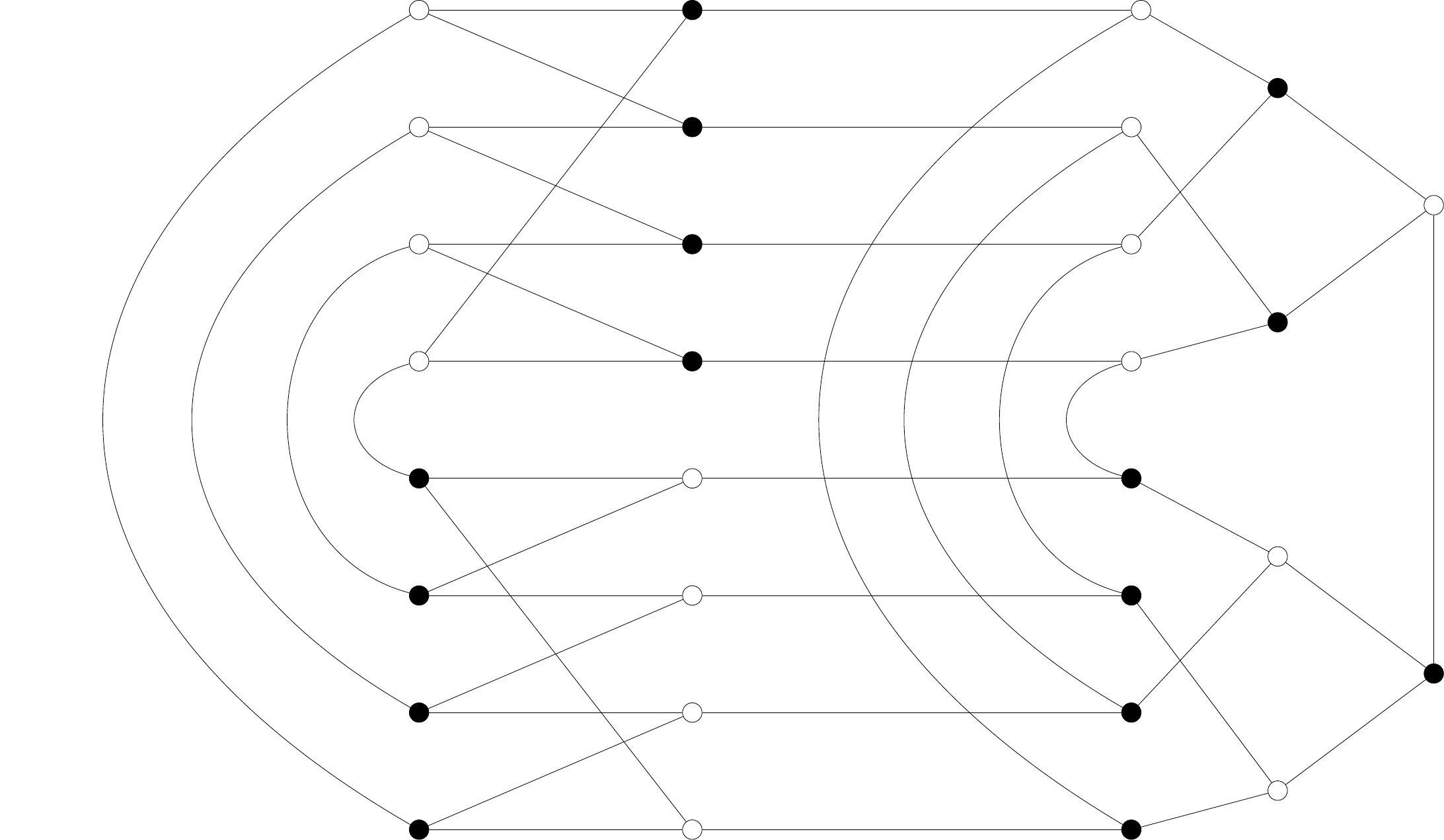} 
    \caption{The counterexample $\mathcal{G}$ on 30 vertices, as constructed in~\cite{AFLR}.}
    \label{fig:config30}
\end{figure}
Another possible drawing of $\mathcal{G}$ is shown in Figure~\ref{fig:config30}, where the \textit{specific join} between the two aforementioned configurations is highlighted and it also allows to explain its automorphism group. More alternative drawings are shown on \textit{Wolfram MathWorld}~\cite{mathworld}.

It follows from~\cite{JG} that the 30-vertex graph $\mathcal{G}$ is the only counterexample to Conjecture~\ref{cong2} up to at least 40 vertices and that there are no counterexamples of girth at least 8 up to at least 48 vertices. Moreover, as a by-product,  Conjecture~\ref{ConjFJLS} and Conjecture~\ref{ConjHeaOnly} are thus verified up to the same orders and girth. 
It is natural to wonder if there are infinitely many counterexamples or further sporadic counterexamples of girth greater than 6 to Conjecture~\ref{cong1} and Conjecture~\ref{cong2}.

The remainder of this manuscript is organised as follows. In Section~\ref{sec:gray}, we show that the Gray graph is also pseudo $2$-factor isomorphic, giving a further counterexample to Conjecture~\ref{cong2}. 
In Section~\ref{sec:comput}, we present the results of our computer searches. Our exhaustive search on cubic bipartite graphs has now been extended to 42 vertices for girth at least 6 and to 52 vertices for girth at least 8, which confirms Conjecture~\ref{ConjFJLS} and Conjecture~\ref{ConjHeaOnly} up to these orders and girth. Moreover, it shows that there are no smaller counterexamples of girth 8 to Conjecture~\ref{cong2} than the Gray graph. 
Next to that, we also verify that there are no further counterexamples among the known censuses of symmetric graphs. 
Finally, in Section~\ref{sec:conclusion} we conclude with some directions for possible future research.

\section{The Gray graph} \label{sec:gray}

\begin{figure}[h!]
    \centering
    \includegraphics[width=.4\textwidth]{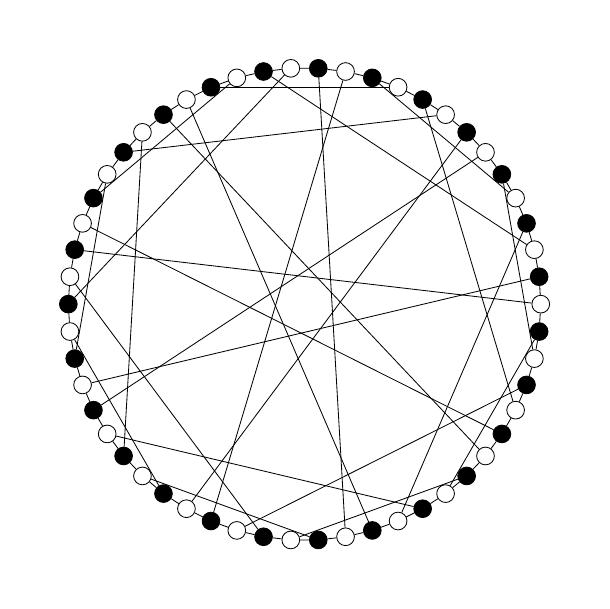} 
    \caption{The Gray graph $\mathcal{GR}$.}
    \label{fig:graygraph}
\end{figure}

In this section, we present another counterexample for Conjecture~\ref{cong2} (and consequently also for Conjecture~\ref{cong1}): \emph{the Gray graph} $\mathcal{GR}$, a cubic bipartite graph on 54 vertices with girth 8, see Figure~\ref{fig:graygraph}. Moreover, it is essentially $4$-edge-connected, has cyclic edge-connectivity 8, and an automorphism group of order 1296.

This graph was originally discovered, but never published, by Marion \linebreak Cameron Gray in 1932. It was re-discovered independently by Bouwer in 1968~\cite{BOU}, in reply to a question posed by  Folkman in 1967~\cite{FOL}.  $\mathcal{GR}$ is interesting as it is the first known example of a cubic graph having the algebraic property of being \textit{semisymmetric}, i.e.\ edge-transitive but not vertex-transitive. In other words, symmetries map every edge to any other edge, but not every vertex to any other vertex; more specifically, vertices in a partition set can only be symmetric to other vertices in the same partition set. In~\cite{MAL} it was shown that the Gray graph $\mathcal{GR}$ is indeed the smallest possible cubic semisymmetric graph.\\

The most common way to construct $\mathcal{GR}$ is described in~\cite{BOU} and consists of the following steps: first, three copies of $K_{3,3}$ are taken, and a chosen edge $e \in E(K_{3,3})$ is subdivided by a vertex in each of the three copies of $K_{3,3}$, and the resulting three vertices are then joined to a new vertex. This is then repeated for each edge in $K_{3,3}$. Other ways to construct $\mathcal{GR}$ may be found in~\cite{MP}, and each one of them explores different and remarkable structural properties of the Gray graph.\\
\begin{figure}[h!]
    \centering
    \includegraphics[width=.9
\textwidth]{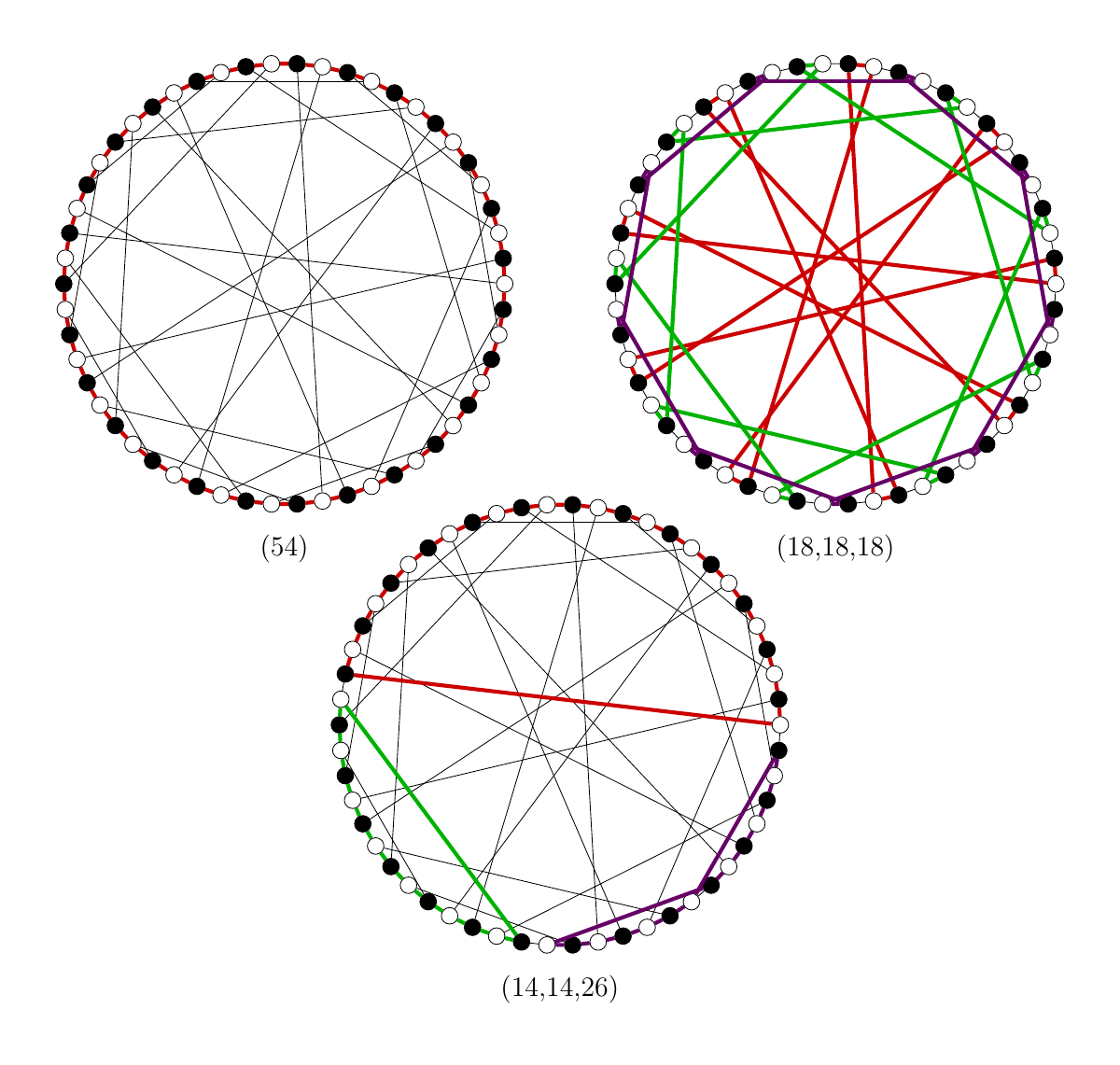} 
    \caption{The Gray graph, with its 2-factor types highlighted.}
    \label{fig:gray2f}
\end{figure}

Using two independent computer programs (see Section~\ref{subsec:computer_details} for details), we determined that there are 10\,752 2-factors in $\mathcal{GR}$, and their cycle sizes are: $(54)$, $(18,18,18)$ and  $(14,14,26)$, as shown in Figure~\ref{fig:gray2f}. These 2-factor types show that the Gray graph is pseudo 2-factor isomorphic but not 2-factor Hamiltonian, so Conjecture~\ref{ConjFJLS} and  Conjecture~\ref{ConjHeaOnly} remain open. The list of all perfect matchings of $\mathcal{GR}$ and the corresponding 2-factors can be found online in~\cite{GitHubRepo}.

\section{Computational results}
\label{sec:comput}
We also extended the computational results from~\cite{JG} in two ways, i.e.\ by exhaustively generating cubic bipartite graphs up to higher orders (cf.\ Section~\ref{subsec:exhaust}) and by exhaustively investigating censuses of symmetrical cubic bipartite graphs (cf.\ Section~\ref{subsec:cencus}). Moreover, in Section~\ref{subsec:computer_details} we give more details on our computer program(s) to test if a given graph is pseudo 2-factor isomorphic and how we verified that our program(s) do not contain any implementation errors.

\subsection{Exhaustive generation}
\label{subsec:exhaust}

In~\cite{JG} the second author generated all cubic bipartite graphs with girth at least 6 up to 40 vertices and all cubic bipartite graphs with girth at least 8 up to 48 vertices and tested which of those graphs is pseudo-2 factor isomorphic. We now extended those searches up to 42 and 52 vertices for girth 6 and 8, respectively. The exact counts are listed in  Table~\ref{table:cubic_bip_graphs}. These counts and the downloadable lists of graphs for smaller orders can also be obtained from \textit{the House of Graphs}~\cite{CDG} at \url{https://houseofgraphs.org/meta-directory/cubic}.

\begin{table}[ht!]
\centering
\begin{tabular}{| c || r | r |}
\hline 
Order & Girth at least 6 & Girth at least 8\\
\hline 
14  &  1  &  \\
16  &  1  &  \\
18  &  3  &  \\
20  &  10  &  \\
22  &  28  &  \\
24  &  162  &  \\
26  &  1 201  &  \\
28  &  11 415  &  \\
30  &  125 571  &  1\\
32  &  1 514 489  &  0\\
34  &  19 503 476  &  1\\
36  &  265 448 847  &  3\\
38  &  3 799 509 760  &  10\\
40  &  57 039 155 060  &  101\\
42  &  \textbf{896 293 917 129}  &  2 510\\
44  &  ?  &  79 605\\
46  &  ?  &  2 607 595\\
48  &  ?  &  81 716 416\\
50  &  ?  &  \textbf{2 472 710 752}\\
52  &  ?  &  \textbf{72 890 068 412}\\
\hline
\end{tabular}

\caption{Counts of all cubic bipartite graphs with girth at least 6 or girth at least~8 for a given order. The counts which are new compared to~\cite{JG} are marked in bold.}

\label{table:cubic_bip_graphs}
\end{table}

The cubic bipartite graphs of girth at least 6 were generated using the generator \texttt{minibaum}~\cite{B}. The generation up to order 42 took approximately 12 CPU years.

The cubic bipartite graphs of girth at least 8 were generated using the generator \texttt{genreg}~\cite{M}, which is a generator for regular graphs. The original version of \texttt{genreg} did not have specific support for generating bipartite graphs, but in the context of the paper~\cite{BGM} Brinkmann recently extended \texttt{genreg} so it could also generate biparitite regular graphs efficiently. In general, \texttt{genreg} is slower than \texttt{minibaum} for generating cubic graphs but it turned out to be significantly faster for generating cubic bipartite graphs of girth at least 8. The generation up to order 52 took approximately 20 CPU years. 

For each of the generated graphs we tested if they are pseudo 2-factor isomorphic, but this did not yield any new (essentially 4-edge-connected) counterexamples to Conjecture~\ref{cong2}. 
(Note that in fact we used \texttt{minibaum}~\cite{B} (and \texttt{genreg}~\cite{M}) to generate cubic bipartite graphs of any connectivity and tested for each of the generated graphs if they are pseudo 2-factor isomorphic or not. This yielded several pseudo 2-factor isomorphic graphs (e.g.\ the known graphs and pseudo 2-factor isomorphic graphs which were obtained by repeated star products). Finally, we tested which of the latter graphs are also essentially 4-edge-connected, which did not yield any new counterexamples.)

In Section~\ref{subsec:computer_details} we give more details about our computer program for testing if a graph is pseudo 2-factor isomorphic. The cost for the pseudo 2-factor isomorphic test was negligible compared to the cost of generating the graphs with \texttt{minibaum} or \texttt{genreg}.
As it is known that any counterexample to Conjecture~\ref{cong2} must have girth at least 6, this leads to the following observations.

\begin{obs}
    The 30-vertex graph $\mathcal{G}$ from Figure~\ref{fig:jan}  is the only counterexample to Conjecture~\ref{cong2} up to at least 42 vertices.
\end{obs}

\begin{obs}
    There are no counterexamples of girth at least 8 to Conjecture~\ref{cong2} up to at least 52 vertices.
\end{obs}

As the Gray graph has 54 vertices and girth 8, this gives us the following Corollary.

\begin{cor}
    The Gray graph is a smallest counterexample of girth at least 8 to Conjecture~\ref{cong2}.
\end{cor}

Since all 2-factor Hamiltonian graphs are pseudo 2-factor isomorphic and $\mathcal{G}$ is not 2-factor Hamiltonian, the above observations also imply the following.

\begin{cor}
Conjecture~\ref{ConjFJLS} and Conjecture~\ref{ConjHeaOnly} hold up to at least 42 vertices and hold
for cubic bipartite graphs with girth at least 8 up to at least 52 vertices.
\end{cor}

\subsection{Investigation of symmetrical graph censuses}
\label{subsec:cencus}

As previously mentioned, the Gray graph is the smallest cubic semisymmetric graph. This naturally leads to the question whether other semisymmetric graphs are pseudo 2-factor isomorphic as well. Conder, Malni\v{c}, Maru\v{s}i\v{c}, and Poto\v{c}nik developed a census of cubic semisymmetric graphs on up to 768 vertices~\cite{Census}, which was later extended on up to 10 000 vertices by Conder and Poto\v{c}nik~\cite{CensusExtended}. In addition, the website \href{https://graphsym.net/}{https://graphsym.net/} contains censuses of highly symmetrical cubic graphs, which can be interesting to check since $K_{3,3}$, the Heawood graph and the Pappus graph are arc-transitive. 
A graph is called \textit{arc-transitive} if its automorphism group acts transitively on its arcs (which are its ordered pairs of adjacent vertices).

Using a computer program (see Section~\ref{subsec:computer_details} for details), we tested whether these censuses of highly symmetrical cubic graphs contain pseudo 2-factor isomorphic graphs different from the known examples. 
We did not find any new pseudo 2-factor isomorphic, nor any new 2-factor Hamiltonian graphs as outlined in \cref{obs:checkP2FI} and~\ref{obs:check2FH}. Recall that a graph $G$ is a \emph{Cayley graph} if there exists a group $\Gamma$ and an inverse-closed subset $S \subseteq \Gamma \setminus \{1_\Gamma\}$
such that the vertices of $G$ are the elements of $\Gamma$, and two vertices
$g,h \in \Gamma$ are adjacent whenever $h = gs$ for some $s \in S$. In this case we denote $G=Cay(\Gamma,S)$.

\begin{obs}\label{obs:checkP2FI}
	There is no essentially 4-edge-connected pseudo 2-factor isomorphic cubic bipartite graph of girth at least 6 different from the Heawood graph, the Pappus graph and the Gray graph which is 
	\begin{itemize}[noitemsep]
		\item vertex-transitive with at most 1 280 vertices,
		\item edge-transitive with at most 10 000 vertices or
		\item Cayley with at most 5 000 vertices.
	\end{itemize}
\end{obs}

Since every 2-factor Hamiltonian graph is also pseudo 2-factor isomorphic, we obtain the following corollary and thus obtain no further counterexamples to Conjecture~\ref{ConjHeaOnly} or Conjecture~\ref{cong2} for these
specific graph families and orders.

\begin{cor}\label{obs:check2FH}
	There is no cyclically 4-edge-connected 2-factor Hamiltonian cubic bipartite graph of girth at least 6 different from the Heawood graph which is
	\begin{itemize}[noitemsep]
		\item vertex-transitive with at most 1 280 vertices,
		\item edge-transitive with at most 10 000 vertices or
		\item Cayley with at most 5 000 vertices.
	\end{itemize}
\end{cor}

\subsection{Computer program details and sanity checks}
\label{subsec:computer_details}

We used two different programs to test if a given graph is pseudo 2-factor isomorphic. To test the exhaustive lists of cubic bipartite graphs of girth at least 6 and girth at least 8 from Section~\ref{subsec:exhaust}, we used the program from~\cite{JG}. This program constructs all perfect matchings and keeps track of the sizes of the corresponding 2-factors (obtained by deleting a perfect matching) and prunes the search if 2-factors with a different parity of the number of cycles have been found (and the graph thus cannot be pseudo 2-factor isomorphic). This straightforward algorithm was fast enough to obtain our computational results from Section~\ref{subsec:exhaust} as the generated graphs are relatively small and the generation process and not the pseudo 2-factor isomorphic test was the bottleneck.

However, in Section~\ref{subsec:cencus} we investigated censuses of highly symmetric graphs for much higher orders (up to 10\,000 vertices). For \emph{non}-pseudo 2-factor isomorphic graphs the program from~\cite{JG} often needed to generate a very large amount of perfect matchings before finding two  where the corresponding 2-factors have a different parity of the number of cycles. Hence, we implemented  a new computer program called \texttt{2FactorParityChecker} -- whose source code can be obtained from~\cite{GitHubRepo} -- that is able to prune much faster on average. The program \texttt{2FactorParityChecker} 
runs multiple threads simultaneously which generate perfect matchings using the randomised Karp-Sipser heuristic method \cite{Kar81} for obtaining an initial -- not necessarily perfect -- matching and extend this matching to a perfect matching using a depth-first search technique.\footnote{We made use of the \texttt{MatchMaker} library~\cite{MatchMaker} for the implementation of the randomised Karp-Sipser and depth-first search method.} 
In addition to these heuristic threads, we have one thread that performs the classic enumeration as implemented in~\cite{JG}. Running these heuristic methods, leads to finding 2-factors with a different parity of the number of cycles much faster.
Thus, a non-pseudo 2-factor isomorphic graph can be processed far more quickly. For graphs that are pseudo 2-factor isomorphic, \texttt{2FactorParityChecker} will generate all perfect matchings with the classic enumeration thread and the heuristic threads will not succeed in finding 2-factors with a different parity of the number of cycles.

Furthermore, as for any computational result, it is important to perform sanity checks to verify the correctness of the implementation of the algorithm. First of all, our implementation can be found on GitHub~\cite{GitHubRepo} as open source software, which can be verified and used by other researchers. Secondly, we compared the output from \texttt{2FactorParityChecker} against the output from the program of~\cite{JG}. We obtained the same perfect matchings and 2-factors for $K_{3,3}$, the Heawood graph, the Pappus graph, Goedgebeur's graph $\mathcal{G}$ and the Gray graph. Moreover, we also adapted both implementations such that they enumerate all perfect matchings -- without terminating early when a 2-factor of different parity is found -- for all graphs of order 30 and girth at least 6 (which are 125\,571 graphs) and obtained the same amount of perfect matchings for each graph.

\section{Concluding remarks}
\label{sec:conclusion}
The first counterexample $\mathcal{G}$ found in 2015 seems hard to generalise into an infinite family of pseudo 2-factor isomorphic graphs; indeed the authors of~\cite{AFLR} pointed out that joining Levi graphs of $n_3$ configurations do not preserve the  property of being pseudo $2$-factor isomorphic. At this point, having found a new counterexample which is a lot more symmetric, we might wonder about the existence of an infinite family of essentially 4-edge-connected cubic bipartite pseudo 2-factor isomorphic graphs. 

We have tried to generalise both the graph $\mathcal{G}$ and the Gray graph without success. One particular approach we tried was through voltage graphs over a group whose regular lifts (or derived graphs) happened to give rise to the special graphs treated in this paper. For definitions of voltage graphs and their regular lifts, please refer to \cite{EJ}. In that context, the theta (multi)graph on two vertices and three edges lifts to $K_{3,3}$ with voltages in $\mathbb{Z}_3$, to the Heawood graph with voltages in $\mathbb{Z}_7$, to the Pappus graph with voltages in $\mathbb{Z}_3^2$, and to the Gray graph with voltages in the semidirect product $\mathbb{Z}_9 \ltimes \mathbb{Z}_3$ (non-abelian group of order 27). 

We made two independent implementations of the algorithm described in~\cite{EJJ} for constructing all regular lifts of girth at least $g$ of a given base graph with voltages in a given group.\footnote{These implementations can also be found on GitHub~\cite{GitHubRepo}.} By executing these algorithms, we were able to conclude that the graph $\mathcal{G}$ is not a regular lift of the theta graph. Moreover, the Pappus graph can be obtained as a regular lift of $K_{3,3}$ with voltage assignments in $\mathbb{Z}_3$, and the Gray graph can also be obtained as a regular lift of the Pappus graph with voltage assignments in $\mathbb{Z}_3$. So we tried to assign voltages to the Gray graph in $\mathbb{Z}_3$ using the aforementioned algorithm, but all of the resulting regular lifts on 162 vertices and girth at least 10 (in fact there was precisely one such graph and it had girth 12) turned out to be not pseudo 2-factor isomorphic. 

This naturally leads to the following problem.
\begin{prob}
    Are there any (essentially 4-edge-connected) cubic bipartite pseudo $2$-factor isomorphic graphs of girth 10 or higher? 
\end{prob}

Please note that Conjecture~\ref{ConjHeaOnly} on $2$-factor Hamiltonian graphs is still open, rendering the Heawood graph still extremely special in this context.

\subsection*{Acknowledgements}
\noindent  The authors are grateful to Louis Stubbe for discussions about regular lifts. 
Jan Goedgebeur and Tibo Van den Eede are supported by Internal Funds of KU Leuven and a grant of the Research Foundation Flanders (FWO) with grant number G0AGX24N.
Jorik Jooken is supported by a Postdoctoral Fellowship of the Research Foundation Flanders (FWO) with grant number 1222524N. 

Several of the computations for this work were carried out using the supercomputer infrastructure provided by the VSC (Flemish Supercomputer Center), funded by the Research Foundation Flanders (FWO) and the Flemish Government. 

This work was also supported by the Italian ``National
Group for Algebraic and Geometric Structures, and their Applications'' (GNSAGA - INdAM).

\end{document}